\documentclass[11pt,reqno,a4paper]{amsart}
\usepackage{amsfonts,amsmath}
\usepackage{fullpage}
\usepackage{graphicx}

\allowdisplaybreaks

\usepackage{stmaryrd}

\newtheorem{theorem}{Theorem}
\newtheorem{lemma}{Lemma}
\newtheorem{coroll}{Corollary}

\DeclareMathOperator{\NegBin}{NegBin}
\DeclareMathOperator{\Geom}{Geom}
\DeclareMathOperator{\Exp}{Exp}

\title{Ordered increasing $k$-trees: introduction and analysis of a preferential attachment network model}

\author[A.~Panholzer]{Alois Panholzer}
\address{Alois Panholzer\\
Institut f{\"u}r Diskrete Mathematik und Geometrie\\
Technische Universit\"at Wien\\
Wiedner Hauptstr. 8-10/104\\
A-1040 Wien, Austria}
\email{Alois.Panholzer@tuwien.ac.at}

\author[G.~Seitz]{Georg Seitz}
\address{Georg Seitz\\
Institut f{\"u}r Diskrete Mathematik und Geometrie\\
Technische Universit\"at Wien\\
Wiedner Hauptstr. 8-10/104\\
A-1040 Wien, Austria}
\email{Georg.Seitz@tuwien.ac.at}

\thanks{This work has been supported by the Austrian Science Foundation FWF, grant S9608-N23.}

\date{\today}


\begin{document}

%
%

\begin{abstract}
   We introduce a random graph model based on $k$-trees, which can be generated by applying a probabilistic
   preferential attachment rule, but which also has a simple combinatorial description.
   We carry out a precise distributional analysis of important parameters for the network model
   such as the degree, the local clustering coefficient and the number of descendants of the nodes and root-to-node distances.
   We do not only obtain results for random nodes, but in particular we also get a precise description of the
   behaviour of parameters for the $j$-th inserted node in a random $k$-tree of size $n$, where $j=j(n)$ might
   grow with $n$. The approach presented is not restricted to this specific $k$-tree model, 
   but can also be applied to other evolving $k$-tree models. 
\end{abstract}

\maketitle


\section{Introduction}

Since the pioneering work of 
\cite{WatStr1998} on real-world networks (as social networks, biological networks or computer networks), various random network models have been introduced that capture at least
part of the typical properties observed frequently. Such properties (see, e.g., \cite{WatStr1998,BarAlb1999}) are, e.g.,
a small average node-to-node distance, a high clustering coefficient, and a power-law degree distribution.

One of the most famous of such random graph models has been introduced by \cite{BarAlb1999}. It uses the idea of 
``preferential attachment'' (or ``success breeds success''), where, starting with a set of nodes, successively nodes are added and linked to a set of nodes by using a specific stochastic growth rule, namely that the probability that a new node is attached to an already existing node is proportional to the degree of that node. A mathematically rigorous definition of this model together with an analysis of important parameters has been given in \cite{BolRio2003}.
It has been pointed out in that work that plane-oriented recursive trees, an important and heavily analyzed tree model
(see, e.g., \cite{MahSmy1995} for a definition and early results), are a special instance of the B\'{a}rabasi-Albert 
graph model.

In this work we introduce a random graph model, which is based on so-called 
$k$-trees\footnote{Here $k \ge 1$ is always an integer. The term $k$-trees, also called $k$-dimensional trees, is somewhat misleading, since, for $k \ge 2$, these graphs are no more trees. 
In particular they should not be confused with $k$-ary trees, which are indeed trees.} (see, e.g., \cite{BeiPip1969,Moo1969}), but where we apply a preferential attachment rule in order to generate them.
Starting with a $k$-clique (a complete connected graph with $k$ vertices) of nodes (the so-called root-clique)
labelled by $0_{1}, 0_{2}, \dots, 0_{k}$, successively the nodes with labels $1, 2, \dots, n$ are inserted, where in each step the inserted node will be attached to all of the nodes of an already existing $k$-clique. But instead of choosing a clique at random we use the probabilistic growth rule that the probability that a new node is attached to an already existing $k$-clique is proportional to one plus the number of nodes that have been previously attached to this $k$-clique (the so-called children of the $k$-clique). In order to also obtain a combinatorial description of these graph families we will consider increasingly labelled ordered $k$-trees and speak about the model of ``ordered increasing $k$-trees''; a precise definition will be given in Section~\ref{sec2}. 

From the construction of $k$-trees it is apparent that for $k=1$ one obtains the model of plane-oriented recursive trees; thus the here studied ordered increasing $k$-trees can be considered as graph families that are generalizations of plane-oriented recursive trees. Quite recently $k$-trees have been introduced as network models in 
\cite{Gao2009,DarHwaBodSor2010} and an analysis of important parameters has been given. In contrast to the 
model we are introducing the considered $k$-trees are there generated by a uniform attachment rule, i.e., 
in each step a new node is attached to a randomly chosen already existing $k$-clique. Combinatorially one might
speak then about the model of ``unordered increasing $k$-trees'', which leads for the special instance $k=1$ 
to the model of (uniform) recursive trees.

We will give a precise distributional analysis of important parameters in ordered increasing $k$-trees 
such as the degree, the local clustering coefficient and the number of descendants of the nodes
and root-to-node distances.
We are here not only interested in a study of quantities for random nodes, but 
a main emphasis is given on describing the behaviour of parameters for the $j$-th inserted node
in a random $k$-tree of size $n$, depending on the growth of $j=j(n)$: we can give a complete characterization
of the limit laws appearing; partially we even obtain exact results. Thus the local behaviour of the nodes during the graph evolution process is described quite well.
Furthermore, using this precise information on the behaviour of the parameters for the $j$-th inserted node we will easily deduce also the limiting behaviour for randomly selected nodes in the $k$-tree (partially we obtain again even exact results). In particular we can show that the distribution of the node-degrees follows asymptotically a so-called power law, i.e., the probability that a randomly selected node has degree $d$ behaves asymptotically as
$\sim c d^{-2-\frac{1}{k}}$, and that the expected local clustering coefficient is rather high 
(e.g., for $k=2$ it is asymptotically, for $n \to \infty$, given by $23 - \frac{9}{4} \pi^{2} \approx 0.793390\dots$).
Moreover the root-to-node distance of node $n$ (but also of a random node), 
is asymptotically Gaussian with expectation $\sim \frac{1}{(k+1)H_{k}} \log n$ and
variance $\sim \frac{H_{k}^{(2)}}{(k+1)H_{k}^{3}} \log n$, where $H_{k} = \sum_{\ell=1}^{n} \frac{1}{\ell}$ and 
$H_{k}^{(2)} = \sum_{\ell=1}^{n} \frac{1}{\ell^{2}}$ denote the first and second order harmonic numbers.

To show our results we use both descriptions of ordered increasing $k$-trees, namely $(i)$ the description via the graph evolution process which often gives rise to a ``bottom-up approach'' when considering the parameter before and after inserting node $n$, and $(ii)$ the combinatorial description as ordered increasing $k$-trees which often
allows a ``top-down approach'' when using a decomposition of the $k$-tree with respect to the root-clique,
see Section~\ref{sec2}. The latter approach
has been applied with success in \cite{DarSor2009,DarHwaBodSor2010} to other $k$-tree models as in particular to
randomly labelled $k$-trees. Both descriptions turn out to be quite useful when analyzing parameters in ordered
increasing $k$-trees; to show our results for the node-degree, the local clustering coefficient and the number of descendants we use the bottom-up approach, whereas for obtaining results on the root-to-node distance we use the top-down approach.

\section{Ordered increasing $k$-trees\label{sec2}}
$k$-trees are families of simple graphs, which have been introduced by \cite{BeiPip1969}.
$k$-trees might be defined recursively in a way analogous to trees: a $k$-tree $T$ is either a $k$-clique (i.e., a complete connected graph with $k$ vertices) or there exists a node $u$ (one might call $u$ endnode), which is incident
to exactly $k$ edges that connect this node to all of the vertices of a $k$-clique, such that, when removing $u$ and the $k$ incident edges from $T$, the remaining graph is itself a $k$-tree. In this paper we will always consider rooted
$k$-trees, which means that in each $k$-tree one $k$-clique is distinguished as the root-clique (the nodes contained
in the root-clique are called root nodes, whereas the remaining nodes are non-root nodes; for the $k$-tree model studied in this work we will also call the non-root nodes ``inserted nodes''). 
Then, apart from the edges connecting the root nodes with each other, this induces a natural orientation on the edges. Thus, for each non-root node,
we can distinguish between ingoing edges (coming from the direction of the root-clique) and outgoing edges, which also defines the in-degree $d^{-}(u)$ and the out-degree $d^{+}(u)$ of a node $u$; for a root node we will only define the out-degree.
It is immediate from the definition that each non-root node $u$ has exactly $k$ ingoing edges, and these edges connect
$u$ with a $k$-clique $K = \{w_{1}, \dots, w_{k}\}$. We might then say that $u$ is a child of the $k$-clique $K$
or that $u$ is attached to $K$ and that $w_{1}, \dots, w_{k}$ are the parents of $u$. For the degree
$d(u)$ of a node $u$ it holds that $d(u) = d^{+}(u) + k$ for a non-root node and $d(u) = d^{+}(u) + k-1$ for a 
root node. We also define the out-degree $d^{+}(K)$ of a $k$-clique $K$ as the number of children of $K$.

Unlike in previously considered $k$-tree models it is for our purpose important to introduce ordered $k$-trees, i.e., we assume that the children of each $k$-clique are linearly ordered (thus one might speak about the first, second, etc.
child of a $k$-clique). Furthermore, we introduce specific labellings of the nodes of ordered $k$-trees, which might
be called increasing labellings (in analogy to the corresponding term for trees, see, e.g., \cite{BerFlaSal1992}).
Given an ordered $k$-tree with $n$ non-root nodes we label the set of root nodes by $\{0_{1}, 0_{2}, \dots, 0_{k}\}$, whereas the non-root nodes are labelled by $\{1, 2, \dots, n\}$ in such a way that the label of a node is always larger than the labels of all its parent nodes (of course, in this context the value of $0_{\ell}$, $1 \le \ell \le k$,
is defined as $0$).
The graph family ``ordered increasing $k$-trees'' can then be described combinatorially as the family of all increasingly labelled ordered $k$-trees. It is apparent from the definition that for $k=1$ one gets the tree family of plane-oriented recursive trees. In what follows we will often use as an abbreviation the term $k$-tree
without further specification, but the meaning should always be ``ordered increasing $k$-tree''.
Furthermore, we will often identify a node with its label, so node $j$ always has the meaning of 
``the node labelled by $j$''.

\begin{figure}
\begin{center}
   \includegraphics[width=12cm]{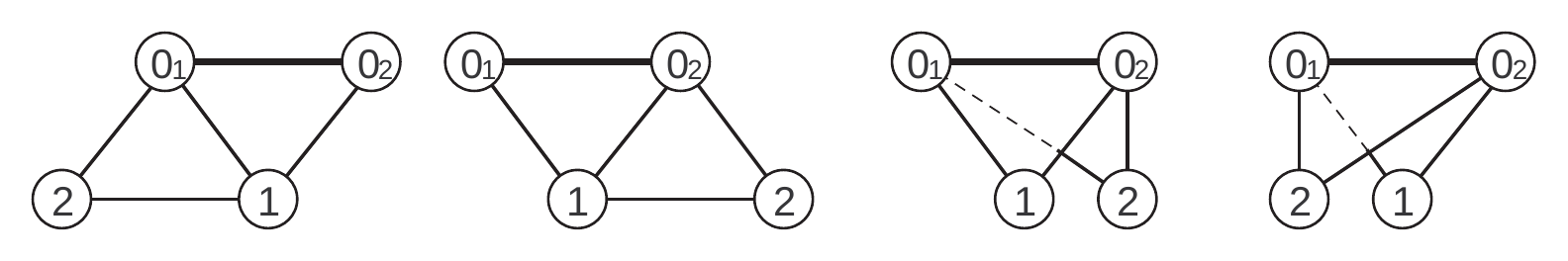}
\end{center}
\vspace*{-5mm}
\caption{All $4$ different $2$-trees of size $2$. In the third and the fourth $2$-tree in the picture the linear order on the children $1$ and $2$ of the root-clique is expressed by drawing $1$ in front of $2$ or vice versa.\label{fig1}}
\end{figure}

Throughout this paper we use the convention that the size $|T|$ of a $k$-tree $T$ is given by the number of non-root nodes; thus the $k$-tree consisting only of the root-clique $K_{0} = \{0_{1}, \dots, 0_{k}\}$ 
has size $0$. Let $T_{n}$ denote the
number of ordered increasing $k$-trees of size $n$ (we do not explicitly express the dependence on $k$, which is of course given). Obviously it holds $T_{0} = T_{1} = 1$. To get an enumeration formula for $T_{n}$ we observe that when inserting a node into a $k$-tree this always increases the number of possible ways of attaching a further node by $k+1$ ($k$ due to the newly generated $k$-cliques and a further one due to a new available position at the
parent $k$-clique). Thus there are always $1+(k+1)(n-1)$ possible ways of inserting node $n$ into a $k$-tree of size
$n-1$. Since each $k$-tree of size $n$ is uniquely obtained from a $k$-tree of size $n-1$ and inserting node $n$ in a possible way it holds that $T_{n} = (1+(k+1)(n-1)) T_{n-1}$, which shows that the number of different ordered increasing $k$-trees of size $n$ is given by
\begin{equation}
   T_{n} = \prod_{\ell=0}^{n-1} \big(1+(k+1) \ell\big) = n! (k+1)^{n} \binom{n-\frac{k}{k+1}}{n}, 
   \quad \text{for $n \ge 0$.}
   \label{eqns2_1}
\end{equation}
In Figure~\ref{fig1} we give all $4$ different $2$-trees of size $2$.

When studying parameters in ordered increasing $k$-trees we always assume the 
``random ordered increasing $k$-tree model'', which means that we assume that each of the $T_{n}$ ordered increasing $k$-trees of size $n$ appears with the same probability.
It is easily seen that for this model the $k$-trees can be obtained indeed by the probabilistic preferential-attachment
growth rule figured out in the introduction. One just has to take in mind that when a $k$-clique $K$ in a $k$-tree
has $\ell$ children, i.e., $d^{+}(K) = \ell$, then there are always exactly $\ell+1$ possible ways of attaching a new node to $K$, namely as the first child, second child, \dots, $(\ell+1)$-th child.
Thus the following evolution process generates ordered increasing $k$-trees uniformly at random:
\begin{itemize}
   \item Step $0$: start with the root clique labelled by $0_{1}, 0_{2}, \dots, 0_{k}$.
   \item Step $n$: the node with label $n$ is attached to any $k$-clique $K$ in the already grown $k$-tree 
   of size $n-1$ with a probability $p(K)$ given by 
   \begin{equation*}
      p(K) = \frac{d^{+}(K)+1}{1+(k+1)(n-1)}.
   \end{equation*}
\end{itemize}

We will also use the combinatorial decomposition of ordered increasing $k$-trees with respect to the root-clique.
To describe this decomposition it is advantagous to introduce two families $\mathcal{T}$ and $\mathcal{S}$ 
of combinatorial objects (they depend on $k$, but we do not explicitly express this). $\mathcal{T}$ is just the family of ordered increasing $k$-trees, whereas $\mathcal{S}$ consists of all ordered increasing $k$-trees, where the root-clique has exactly one child.
Of course, an object of $\mathcal{T}$, where the root-clique has exactly $\ell$ children, can be obtained, after identification of the root nodes and an order-preserving relabelling, by a sequence of $\ell$ objects of $\mathcal{S}$.
Furthermore, when considering objects in $\mathcal{S}$ the child-node attached to the root-clique has to be labelled by $1$, and together with all choices of $k-1$ nodes from the root-clique it is forming exactly $k$ different $k$-cliques, which, after relabelling, can themselves be considered as root-cliques of objects of $\mathcal{T}$.
Thus we obtain the following formal description of the families $\mathcal{T}$ and $\mathcal{S}$
(see, e.g., \cite{FlaSed2009} for an explanation of such formal specifications):
\begin{equation}
   \label{eqns2_3}
   \begin{split}
   \mathcal{T} & = \mathcal{S}^{0} \; \dot{\cup} \; \mathcal{S}^{1} \; \dot{\cup} \; \mathcal{S}^{2} \; 
   \dot{\cup} \; \mathcal{S}^{3} \; \dot{\cup} \; \cdots \; = \; \text{\textsc{Seq}}(\mathcal{S}), \\
   \mathcal{S} & = \{1\} \times \mathcal{T}^{k} = \mathcal{Z}^{\boxempty} \ast \mathcal{T}^{k}.
   \end{split}
\end{equation}
When denoting by $T_{n}$ and $S_{n}$ the number of objects in the families $\mathcal{T}$ and $\mathcal{S}$,
respecively, of size $n$ and by $T(z) := \sum_{n \ge 0} T_{n} \frac{z^{n}}{n!}$ and $S(z) := \sum_{n \ge  0}
S_{n} \frac{z^{n}}{n!}$ their exponential generating functions, we obtain by using the symbolic method
(see again, e.g., \cite{FlaSed2009}) immediately the following system of equations:
\begin{equation*}
   T(z) = \frac{1}{1-S(z)}, \quad S'(z) = T(z)^{k}, \quad S(0)=0,
\end{equation*}
which has the solution
\begin{equation}
   T(z) = \big(1-(k+1)z\big)^{-\frac{1}{k+1}} \quad \text{and} \quad
   S(z) = 1-\big(1-(k+1)z\big)^{\frac{1}{k+1}}.
   \label{eqns2_4}
\end{equation}
Extracting coefficients from $T(z)$ and $S(z)$ shows again that the number $T_{n}$ of ordered increasing $k$-trees
of size $n$ is given by \eqref{eqns2_1}, whereas $S_{n} = (n-1)! (k+1)^{n-1} \binom{n-1-\frac{1}{k+1}}{n-1}$, 
$n \ge 1$.

\section{Parameters studied and results}

\subsection{Parameters studied}
Next we give a definition of the quantities studied in random ordered increasing $k$-trees.
For better readability we do not explicitly express the dependence of the quantities on $k$, which is of course always given, in the notations.

The r.v. $Y_{n,j}$ counts the out-degree (see Section~\ref{sec2}) of node $j$ in a random $k$-tree of size $n$,
whereas the r.v. $\bar{Y}_{n}$ counts the out-degree of a random inserted node in a random $k$-tree of size $n$,
i.e., the out-degree of a node picked at random from the non-root nodes $\{1, 2, \dots, n\}$.
The r.v. $Y_{n,0}$ counts the out-degree of the root-node $0_{1}$ in a random $k$-tree of size $n$; of course, due to symmetry, the corresponding r.v. are identically distributed for each of the root-nodes $0_{1}, \dots, 0_{k}$
and do not have to be introduced separately.

The r.v. $C_{n,j}$ counts the local clustering coefficient of node $j$ in a random $k$-tree of size $n$.
The local clustering coefficient has been introduced by \cite{WatStr1998} and is considered 
as an important parameter in the study of real-world networks.
The local clustering coefficient $C_{G}(u)$ of a node $u$ in a graph $G(V,E)$ is defined as the proportion of
edges between neighbours of $u$ divided by the number of edges between the neighbours that could possibly exist;
formally $C_{G}(u)$ is given by
\begin{equation}
   C_{G}(u) = \begin{cases}
             \frac{|\{e \in E \: : \: e=(x,y) \; \text{with} \; x,y \in N(u)\}|}{\binom{d(u)}{2}}, 
             \quad \text{if $d(u) \ge 2$,} \\
             0, \quad \text{if $d(u)=0$ or $d(u)=1$,}
          \end{cases}
   \label{eqns3_1}
\end{equation}
where $N(u)$ denotes the set of neighbours (i.e., adjacent nodes) of $u$.
The r.v. $\bar{C}_{n}$ counts the local clustering coefficient of a randomly selected node 
(amongst the root nodes and the inserted nodes) in a random $k$-tree of size $n$.

The r.v. $X_{n,j}$ counts the number of descendants of node $j$ in a random $k$-tree of size $n$.
Whether a node $w$ is a descendant of $u$ might be defined recursively: $w$ is a descendant of $u$ if either $u=w$
or if $w$ has a parent node, which is a descendant of $u$. One might then also say that $u$ is an ancestor of $w$.
The r.v. $\bar{X}_{n}$ counts the number of descendants of a random inserted node in a random $k$-tree of size $n$.

The r.v. $D_{n}$ counts the distance between the root node $0_{1}$ and node $n$ in a random $k$-tree of size $n$.
As usual the distance between nodes in a graph is measured by the minimal number of edges contained in a path
amongst all paths connecting these nodes. Again the corresponding r.v. are identically distributed for each of the
root-nodes $0_{1}, \dots, 0_{k}$.
Furthermore, as a direct consequence of the evolution process of random $k$-trees one obtains that the distance
$D_{n,j}$ between the root node $0_{1}$ and node $j$ in a random $k$-tree of size $n$ is distributed as $D_{j}$
and thus does not have to be studied separately.
The r.v. $\bar{D}_{n}$ counts the distance between the root node $0_{1}$ and a random inserted node in a random $k$-tree of size $n$.

\begin{figure}
\begin{center}
   \includegraphics[width=6cm]{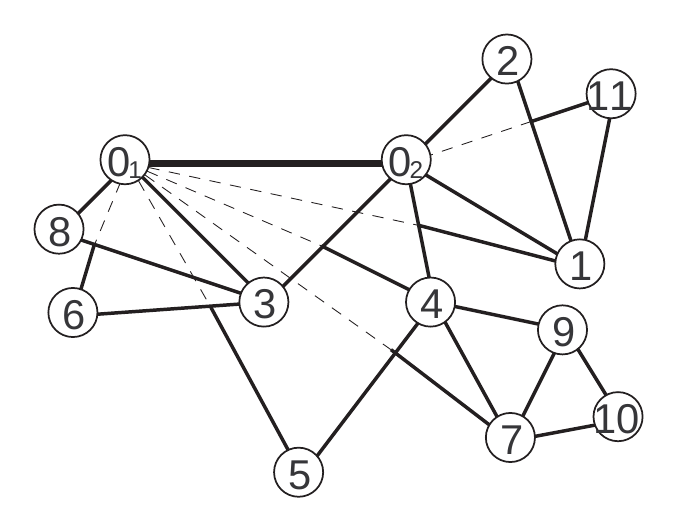}
\end{center}
\vspace*{-8mm}
\caption{An example of a $2$-tree of size $11$. Node $4$ has out-degree three and five descendants (counting the node as a descendant of itself). The local clustering coefficient of node $4$ is $0.4$, since there are four edges between the five neighbours of this node (see the definition of the local clustering coefficient). 
The distance of node $4$ to the root node $0_{1}$ is one.\label{fig2}}
\end{figure}

In Figure~\ref{fig2} we give an example of a $2$-tree together with the parameters studied.

\subsection{Results}

\subsubsection*{Degree of the nodes}

\begin{theorem}
   \label{the1}
   The r.v. $Y_{n,j}$, which counts the out-degree of node $j$ in a random $k$-tree of size $n$, 
   has the following exact distribution:
   \begin{equation*}
      \mathbb{P}\{Y_{n,j} = m\} = \frac{\binom{j-\frac{k}{k+1}}{j}}{\binom{n-\frac{k}{k+1}}{n} \binom{n}{j}}
      \sum_{\ell=0}^{m} \binom{m}{\ell} (-1)^{\ell} \binom{n-\frac{k(2+\ell)}{k+1}}{n-j}, \quad 
      \text{for $n \ge j \ge 1$ and $m \ge 0$.}
   \end{equation*}
   
   The limiting distribution behaviour of $Y_{n,j}$ is, for $n \to \infty$ and depending on the growth of $j$,
   characterized as follows.
\begin{itemize}
\item The region for $j$ fixed.
The normalized random variable $n^{-\frac{k}{k+1}} Y_{n,j}$ converges in distribution to a r.v. $Y_{j}$, 
i.e., $n^{-\frac{k}{k+1}} Y_{n,j} \xrightarrow{(d)} Y_{j}$, which is fully characterized by its moments.
The $s$-th moments of $Y_{j}$ are, for $s \ge 0$, given by
\begin{equation*}
   \mathbb{E}(Y_{j}^{s}) = \frac{s! \Gamma(j+\frac{1}{k+1})}{\Gamma(j+1+\frac{k(s-1)}{k+1})}.
\end{equation*}

\item The region for $j$ small: $j \to \infty$ such that $j = o(n)$.
The normalized random variable $\big(\frac{j}{n}\big)^{\frac{k}{k+1}}Y_{n,j}$ 
is asymptotically exponentially distributed with parameter $1$,
$\big(\frac{j}{n}\big)^{\frac{k}{k+1}} Y_{n,j} \xrightarrow{(d)} \Exp(1)$, i.e., 
$\big(\frac{j}{n}\big)^{\frac{k}{k+1}} Y_{n,j} \xrightarrow{(d)} Y$,
where the $s$-th moments of $Y$ are, for $s \ge 0$, given by
\begin{equation*}
   \mathbb{E}(Y^{s}) = s!.
\end{equation*}

\item The central region for $j$: $j \to \infty$ such that $j \sim \rho n$, with $0 < \rho < 1$.
The random variable $Y_{n,j}$ is asymptotically geometrically distributed with success probability $\rho^{\frac{k}{k+1}}$, 
$Y_{n,j} \xrightarrow{(d)} \Geom(\rho^{\frac{k}{k+1}})$, i.e., $Y_{n,j} \xrightarrow{(d)} Y_{\rho}$, 
where the probability mass function of $Y_{\rho}$ is given by
\begin{equation*}
   \mathbb{P}\{Y_{\rho} = m\} = \rho^{\frac{k}{k+1}} \big(1-\rho^{\frac{k}{k+1}}\big)^{m},
   \quad \text{for} \; m \ge 0.
\end{equation*}

\item The region for $j$ large: $j \to \infty$ such that $\tilde{j} := n-j = o(n)$.
It holds that $\mathbb{P}\{Y_{n,j} = 0\} \to 1$.
\end{itemize}
\end{theorem}

\begin{theorem}
   \label{the2}
   The r.v. $Y_{n,0}$, which counts the out-degree of the root node $0_{1}$ in a random $k$-tree of
   size $n$, has the following exact distribution:
   \begin{equation*}
      \mathbb{P}\{\bar{Y}_{n,0} = m\} = \frac{\binom{m-\frac{k-1}{k}}{m}}{\binom{n-\frac{k}{k+1}}{n}}
      \sum_{\ell=0}^{m} \binom{m}{\ell} (-1)^{\ell} \binom{n-1-\frac{k\ell}{k+1}}{n}, \quad \text{for $m \ge 0$.}
   \end{equation*}
   For $n \to \infty$, the normalized random variable $n^{-\frac{k}{k+1}} Y_{n,0}$ converges in distribution to a 
   r.v. $Y_{0}$, i.e., $n^{-\frac{k}{k+1}} Y_{n,0} \xrightarrow{(d)} Y_{0}$, 
   which is fully characterized by its moments.
   The $s$-th moments of $Y_{0}$ are, for $s \ge 0$, given by
   \begin{equation*}
      \mathbb{E}(Y_{0}^{s}) = \frac{\Gamma(\frac{1}{k+1}) 
      \Gamma(s+\frac{1}{k})}{\Gamma(\frac{1}{k}) \Gamma(\frac{k}{k+1}s+\frac{1}{k+1})}.
   \end{equation*}
\end{theorem}

\begin{theorem}
   \label{the3}
   The r.v. $\bar{Y}_{n}$, which counts the out-degree of a random inserted node in a random $k$-tree of
   size $n$, has the following exact distribution:
   \begin{equation*}
      \mathbb{P}\{\bar{Y}_{n} = m\} = \frac{1}{n \binom{n-\frac{k}{k+1}}{n}} \sum_{\ell=0}^{m}
      \frac{\binom{m}{\ell} (-1)^{\ell}}{k(\ell+2)+1}
      \left(\binom{n+\frac{1}{k+1}}{n}-\binom{n-\frac{k(\ell+2)}{k+1}}{n}\right), \quad \text{for $m \ge 0$.}
   \end{equation*}
   For $n \to \infty$, $\bar{Y}_{n}$ converges in distribution to a discrete r.v. $\bar{Y}$, i.e., 
   $\bar{Y}_{n} \xrightarrow{(d)} \bar{Y}$, with
   \begin{equation*}
      \mathbb{P}\{\bar{Y}=m\} = p_{m} := \frac{k+1}{k (m+1) \binom{m+2+\frac{1}{k}}{m+1}}, \quad \text{for $m \ge 0$.}
   \end{equation*}
   Since $p_{m} \sim \frac{k+1}{k} \Gamma\big(2+\frac{1}{k}\big) m^{-2-\frac{1}{k}}$, for $m \to \infty$,
   it follows that $\bar{Y}_{n}$ follows asymptotically a power-law distribution with exponent $2+\frac{1}{k}$.
\end{theorem}

\subsubsection*{Local clustering coefficient}

\begin{lemma}
   \label{lem1}
   For any $k$-tree $T$ the local clustering coefficient $C_{T}(u)$ of a node $u$ only depends on the 
   degree $d(u)$ of $u$. It holds then for $d(u) \ge k \ge 2$:
   \begin{equation*}
      C_{T}(u) = \frac{2(k-1)}{d(u)} - \frac{(k-1)(k-2)}{d(u) (d(u)-1)}.
   \end{equation*}
\end{lemma}

\begin{theorem}
   \label{the4}
   Let the r.v. $\bar{C}_{n}$ count the local clustering coefficient of a random node in a random $k$-tree of
   size $n$. Then the expected local clustering coefficient $\mathbb{E}(\bar{C}_{n})$ behaves, for $n \to \infty$,
   as follows (here $\Psi(x) = (\ln \Gamma(x))'$ denotes the Psi-function and $\Psi'(x)$ its derivative):
   \begin{align*}
      & \mathbb{E}(\bar{C}_{n}) \to c_{k} := \sum_{m \ge 0} \frac{(k+1)(k-1)}{k(m+1)(m+k)\binom{m+2+\frac{1}{k}}{m+1}}
      \Big(2-\frac{k-2}{m+k-1}\Big) \\
      & = (k+1) \bigg(\frac{1}{k-1} + \frac{3}{k}\binom{k-4-\frac{1}{k}}{k-2}\sum_{\ell=1}^{k-2}
      \frac{1}{\ell^{2} \binom{\ell+3+\frac{1}{k}-k}{\ell}} + \frac{3}{k} \sum_{\ell=0}^{k-3} \frac{(-1)^{\ell}}
      {(\ell+1)(\ell-1-\frac{1}{k}) \binom{k-3}{\ell}} \\
      & \qquad \qquad \quad \mbox{} - \frac{3}{k} \binom{k-4-\frac{1}{k}}{k-2} \Psi'\big(4+\frac{1}{k}-k\big)\bigg).
   \end{align*}
   It further holds that $c_{k} \to 1$, for $k \to \infty$.
\end{theorem}

\begin{table}
  \begin{equation*}
  \begin{array}{|c|c|c|}
    \hline
    \rule[-1.5ex]{0pt}{4.5ex} k & c_{k} = \lim\limits_{n \to \infty} \mathbb{E}(\bar{C}_{n}) & \text{numerically} \\
    \hline
    \rule[-1.5ex]{0pt}{4.5ex} 2 & 23 - \frac{9}{4} \pi^{2} & 0.793390\dots \\
    \hline
    \rule[-1.5ex]{0pt}{4.5ex} 3 & -5 + \frac{16}{3} \Psi'(\frac{4}{3}) & 0.843184\dots \\
    \hline
    \rule[-1.5ex]{0pt}{4.5ex} 4 & \frac{1051}{96} - \frac{75}{128} \Psi'(\frac{1}{4}) & 0.871356\dots \\
    \hline
    \rule[-1.5ex]{0pt}{4.5ex} 5 & \frac{512}{125} - \frac{72}{625} \Psi'(-\frac{4}{5}) & 0.889998\dots \\
    \hline
  \end{array}
  \;
  \begin{array}{|c|c|c|}
    \hline
    \rule[-1.5ex]{0pt}{4.5ex} k & c_{k} = \lim\limits_{n \to \infty} \mathbb{E}(\bar{C}_{n}) & \text{numerically} \\
    \hline
    \rule[-1.5ex]{0pt}{4.5ex} 6 & \frac{148003}{57024} - \frac{2695}{62208} \Psi'(-\frac{11}{6}) & 0.903449\dots \\
    \hline
    \rule[-1.5ex]{0pt}{4.5ex} 10 & & 0.933975\dots \\
    \hline
    \rule[-1.5ex]{0pt}{4.5ex} 50 & & 0.982804\dots \\
    \hline
    \rule[-1.5ex]{0pt}{4.5ex} 100 & & 0.990885\dots \\
    \hline
  \end{array}
  \end{equation*}
  \caption{The limit $c_{k}$ of the expected local clustering coefficient $\mathbb{E}(\bar{C}_{n})$ 
  for small values of $k$.}
\end{table}

\subsubsection*{Number of descendants}

\begin{theorem}
   \label{the5}
   The r.v. $X_{n,j}$, which counts the number of descendants of node $j$ in a random $k$-tree of size $n$, 
   has the following exact distribution:
   \begin{equation*}
      \mathbb{P}\{X_{n,j} = m\} = \frac{\binom{m-1-\frac{1}{k+1}}{m-1} \binom{n-m-1+\frac{2}{k+1}}{n-m-j+1}}
      {\binom{n-\frac{k}{k+1}}{n-j}}, \quad \text{for $n \ge j \ge 1$ and $m \ge 1$.}
   \end{equation*}
   
   The limiting distribution behaviour of $X_{n,j}$ is, for $n \to \infty$ and depending on the growth of $j$,
   characterized as follows.
\begin{itemize}
\item The region for $j$ fixed.
The normalized random variable $\frac{X_{n,j}}{n}$ is asymptotically Beta-distributed,
$\frac{X_{n,j}}{n} \xrightarrow{(d)} \beta(\frac{k}{k+1},j-1+\frac{2}{k+1})$, i.e.,
$\frac{X_{n,j}}{n} \xrightarrow{(d)} X_{j}$,
where the $s$-th moments of $X_{j}$ are, for $s \ge 0$, given by
\begin{equation*}
   \mathbb{E}(X_{j}^{s}) = \frac{\binom{s-\frac{1}{k+1}}{s}}{\binom{s+j-\frac{k}{k+1}}{s}}.
\end{equation*}

\item The region for $j$ small: $j \to \infty$ such that $j = o(n)$.
The normalized random variable $\frac{j}{n}X_{n,j}$ is asymptotically Gamma-distributed,
$\frac{j}{n}X_{n,j} \xrightarrow{(d)} \gamma(\frac{k}{k+1},1)$, i.e., $\frac{j}{n}X_{n,j} \xrightarrow{(d)} X$,
where the $s$-th moments of $X$ are, for $s \ge 0$, given by
\begin{equation*}
   \mathbb{E}(X^{s}) = s! \binom{s-\frac{1}{k+1}}{s}.
\end{equation*}

\item The central region for $j$: $j \to \infty$ such that $j \sim \rho n$, with $0 < \rho < 1$.
The shifted random variable $X_{n,j}-1$ is asymptotically negative binomial-distributed, 
$X_{n,j} -1 \xrightarrow{(d)} \NegBin(\frac{k}{k+1},\rho)$, i.e., $X_{n,j} -1 \xrightarrow{(d)} X_{\rho}$, 
where the probability mass function of $X_{\rho}$ is given by
\begin{equation*}
   \mathbb{P}\{X_{\rho} = m\} = \binom{m-\frac{1}{k+1}}{m} \rho^{\frac{k}{k+1}} (1-\rho)^{m},
   \quad \text{for} \; m \ge 0.
\end{equation*}

\item The region for $j$ large: $j \to \infty$ such that $\tilde{j} := n-j = o(n)$.
It holds that $\mathbb{P}\{X_{n,j} = 1\} \to 1$.
\end{itemize}
\end{theorem}

\begin{theorem}
   \label{the6}
   The r.v. $\bar{X}_{n}$, which counts the number of descendants of a random inserted node in a random $k$-tree of
   size $n$, has the following exact distribution (with $m \ge 1$):
   \begin{equation*}
      \mathbb{P}\{\bar{X}_{n} = m\} = \frac{\binom{m-1-\frac{1}{k+1}}{m-1}}{n \binom{n-\frac{k}{k+1}}{n}}
      \sum_{\ell = 0}^{m-1} \frac{\binom{m-1}{\ell} (-1)^{\ell}}{(k+1)(\ell+1)+k}
      \left(\binom{n+\frac{1}{k+1}}{n}-\binom{n-\ell-2+\frac{2}{k+1}}{n}\right).
   \end{equation*}
   For $n \to \infty$, $\bar{X}_{n}$ converges in distribution to a discrete r.v. $\bar{X}$, i.e., 
   $\bar{X}_{n} \xrightarrow{(d)} \bar{X}$, with
   \begin{equation*}
      \mathbb{P}\{\bar{X}=m\} = \frac{k}{(k+1)(m+\frac{k}{k+1})(m-\frac{1}{k+1})}, \quad \text{for $m \ge 1$.}
   \end{equation*}
\end{theorem}

\subsubsection*{Root-to-node distance}

\begin{theorem}
   \label{the7}
   The r.v. $D_{n}$, which measures the distance between the root node $0_{1}$ and node $n$ in a random $k$-tree 
   of size $n$, is, for $n \to \infty$, asymptotically Gaussian distributed, where the rate of convergence 
   is of order $\mathcal{O}\big(\frac{1}{\sqrt{\log n}}\big)$:
   \begin{equation*}
      \sup_{x \in \mathbb{R}} \: \Bigg| \:
      \mathbb{P}\Bigg\{\frac{D_{n} - \mathbb{E}(D_{n})}
		  {\sqrt{\mathbb{V}(D_{n})}} \le x\Bigg\} - \Phi(x) 
		  \: \Bigg|
      = \mathcal{O}\Big(\frac{1}{\sqrt{\log n}}\Big),
   \end{equation*}
   and the expectation $\mathbb{E}(D_{n})$ and the variance 
   $\mathbb{V}(D_{n})$ satisfy
   \begin{equation*}
      \mathbb{E}(D_{n}) = \frac{1}{(k+1) H_{k}} \log n + \mathcal{O}(1), \quad
   	  \mathbb{V}(D_{n}) = \frac{H_{k}^{(2)}}{(k+1) H_{k}^{3}} \log n + \mathcal{O}(1).
   \end{equation*}
\end{theorem}

\begin{coroll}
   \label{cor1}
   The r.v. $\bar{D}_{n}$, which measures the distance between the root node $0_{1}$ and a random inserted node in a
   random $k$-tree of size $n$, is, for $n \to \infty$, asymptotically Gaussian distributed:
   $\mathbb{P}\left\{\frac{\bar{D}_{n} - \mathbb{E}(\bar{D}_{n})}{\sqrt{\mathbb{V}(\bar{D}_{n})}} \le x\right\} 
   \to \Phi(x)$, for all $x \in \mathbb{R}$,
   with expectation and variance satisfying $\mathbb{E}(\bar{D}_{n}) = \frac{1}{(k+1) H_{k}} \log n + \mathcal{O}(1)$
   and $\mathbb{V}(\bar{D}_{n}) = \frac{H_{k}^{(2)}}{(k+1) H_{k}^{3}} \log n + \mathcal{O}(1)$.
\end{coroll}
Here $\Phi(x)$ denotes the distribution function of the standard normal distribution $\mathcal{N}(0,1)$.

\section{Brief outline of the proof of the results}

\subsection{Degree of the nodes\label{ssec41}}
In order to get a suitable description of the r.v. $Y_{n,j}$ we consider the graph evolution process of $k$-trees.
The following observation is crucial to our approach: each node $x$ attached to node $j$ increases the number of
possibilities of attaching a new node to $j$ by exactly $k$ ($1$ possibility more at the $k$-clique where node $x$
is attached and $k-1$ possibilities more due to the $k-1$ new $k$-cliques containing $x$ and $j$).
Thus if node $j \ge 1$ has out-degree $m$ there are exactly $(m+1) k$ possibilities of attaching a new node that
increases the out-degree of node $j$, whereas the remaining possibilities will keep the out-degree unchanged. 

Thus if we count by $T_{n,j,m} := T_{n} \mathbb{P}\{Y_{n,j} = m\}$ the number of $k$-trees of size $n$ such that node
$j$ has out-degree $m$, we immediately get the following recurrence:
\begin{equation*}
   T_{n,j,m} = \big((k+1) n - km - 2k\big) T_{n-1,j,m} + km T_{n-1,j,m-1}, \quad \text{for $n > j \ge 1$ and $m \ge 0$},
\end{equation*}
with $T_{j,j,0} = T_{j}$, for $j \ge 1$, and $T_{j,j,m} = 0$, for $m > 0$.
Introducing the generating functions 
$T^{[j]}(z,v) := \sum_{n \ge j} \sum_{m \ge 0} T_{n,j,m} \frac{z^{n-j}}{(n-j)!} v^{m}$ leads to the following 
linear first order partial differential equation:
\begin{equation*}
   \big(1-(k+1)z\big) T_{z}^{[j]}(z,v) + kv (1-v) T_{v}^{[j]}(z,v) - \big(k(j-1+v)+j+1\big) T^{[j]}(z,v) = 0,
   \quad T^{[j]}(0,v)=T_{j},
\end{equation*}
which can be solved by applying the method of characteristics. The solution is given by the following expression:
\begin{equation*}
   T^{[j]}(z,v) = \frac{T_{j}}{\big(1-v\big(1-(1-(k+1)z)^{\frac{k}{k+1}}\big)\big)
   \big(1-(k+1)z\big)^{\frac{kj-k+j+1}{k+1}}},
\end{equation*}
and extracting coefficients immediately shows the exact formula for the probabilities $\mathbb{P}\{Y_{n,j}=m\}$ 
given in Theorem~\ref{the1}.
To show the limiting distribution results given in Theorem~\ref{the1} we use, depending on the growth behaviour of
$j=j(n)$, different approaches (see \cite{KubPan2007} for similar considerations on the node-degree of 
increasing trees). 
For the two cases $j$ fixed and $j \to \infty$, such that $j = o(n)$, we use the method
of moments, where we study the explicit expression for the $s$-th factorial moments obtained after extracting coefficients from the $s$-th derivative of $T^{[j]}(z,v)$ w.r.t. $v$ evaluated at $v=1$.
For the remaining two cases $j \sim \rho n$, with $0 < \rho < 1$, and $n-j = o(n)$ we directly study the exact
expression for the probabilities. 

To show Theorem~\ref{the2} concerning the out-degree of the root node $0_{1}$ one can use the same approach as for 
a non-root $j$, but one has to start with a slightly different recurrence. The asymptotic considerations are similar to the case $j$ fixed.

For obtaining the results given in Theorem~\ref{the3} one simply uses the relation
$\mathbb{P}\{\bar{Y}_{n} = m\} = \frac{1}{n} \sum_{j=1}^{n} \mathbb{P}\{Y_{n,j} = m\}$ and Theorem~\ref{the1}. 
In order to get the explicit expression for the probabilities given in the theorem we use a
hypergeometric identity for simplifying it.

\subsection{Local clustering coefficient}
The crucial observation for analyzing the local clustering coefficient in $k$-trees is 
that the local clustering coefficient $C_{T}(u)$ of a node $u$ in a $k$-tree $T$ only depends on the degree $d(u)$ 
of the corresponding node; the exact relation is expressed in Lemma~\ref{lem1}.
To show this we will, according to the definition~\eqref{eqns3_1}, count the number $M(u)$ of edges between 
neighbours of $u$. Consider a node $u$ in a $k$-tree; then it always holds that $d(u) \ge k-1$.
If $d(u) = k-1$ then the $k$-tree can consist only of a single root-clique and $u$ is one of the root nodes; 
thus all $k-1$ neighbours of $u$ are connected with each other, which implies $M(u) = \binom{k-1}{2}$.
In order to determine $M(u)$ when $d(u) \ge k$ we observe that in any $k$-tree holds that when increasing the degree of a node $u$ by $1$ then the number of edges between neighbours
of $u$ increases exactly by $k-1$; this holds since a new node $w$ adjacent to $u$ generates a $k$-clique,
such that $w$ is also adjacent to $k-1$ neighbours of $u$.
Thus $M(u) = \binom{k-1}{2} + (k-1)(d(u)-k+1)$, for $d(u) \ge k-1$, which implies Lemma~\ref{lem1}.

Of course, due to Lemma~\ref{lem1}, one can immediately obtain distributional relations between r.v. 
measuring the degree (or out-degree) and the local clustering coefficient of nodes in $k$-trees.
In particular we are interested in the r.v. $\bar{C}_{n}$ measuring the local clustering coefficient of a random
node in a random $k$-tree of size $n$ (of course, similar considerations for the local clustering coefficient of
node $j$ can be made also, but we skip them here). One gets then
\begin{equation*}
   \bar{C}_{n} \stackrel{(d)}{=} \frac{2(k-1)}{\tilde{Y}_{n}} - \frac{(k-1)(k-2)}{\tilde{Y_{n}} (\tilde{Y}_{n}-1)},
\end{equation*}
where $\tilde{Y}_{n}$ measures the degree of a randomly selected node (amongst the root nodes and inserted nodes)
in a $k$-tree of size $n$.
Of course, the distribution of $\tilde{Y}_{n}$, and thus also the distribution of $\bar{C}_{n}$
is fully determined by the previously studied r.v. $\bar{Y}_{n}$ and $Y_{n,0}$. In particular it easily follows that $\tilde{Y}_{n} \xrightarrow{(d)} \bar{Y}+k$, where the distribution of the discrete random variable $\bar{Y}$ 
is characterized in Theorem~\ref{the3}.
The main quantity of interest in this context is the expected local clustering coefficient of a random node.
Since $\tilde{Y}_{n} \xrightarrow{(d)} \bar{Y}+k$, with $\bar{Y}$ a discrete r.v., and since the function
$f(m) = \frac{2(k-1)}{m} - \frac{(k-1)(k-2)}{m(m-1)}$ is uniformly bounded for $m \ge k$, it immediately follows
that
\begin{equation*}
   \mathbb{E}(\bar{C}_{n}) \to c_{k} := \sum_{m \ge k} \mathbb{P}\{\bar{Y}+k=m\} 
   \Big(\frac{2(k-1)}{m}-\frac{(k-1)(k-2)}{m(m-1)}\Big),
\end{equation*}
which leads to the first expression for $c_{k}$ given in Theorem~\ref{the4}. The second one, which is advantageous
when computing $c_{k}$ for small $k$, can be obtained by rather lengthy manipulations with  beta integrals
and their derivatives.

\subsection{Number of descendants}
For a recursive description of the r.v. $X_{n,j}$ we consider the graph evolution process of $k$-trees.
Here the following observation is crucial: each node $x$ attached to a descendant of node $j$ increases the number of
possibilities of attaching a new node to a descendant of $j$ by exactly $k+1$ ($1$ possibility more at the 
$k$-clique where node $x$ is attached and $k$ possibilities more due to the $k$ new $k$-cliques containing $x$).
Thus if node $j \ge 1$ has $m$ descendants there are exactly $(k+1)m-1$ possibilities of attaching a new node that
increases the number of descendants of node $j$, whereas the remaining possibilities will keep the number of descendants unchanged. 

This description allows a recursive approach analogeous to the one sketched in Subsection~\ref{ssec41}. Also the asymptotic considerations are very similar to the ones discussed there, but somewhat simpler, due to the 
closed formul{\ae} for the exact results.

\subsection{Root-to-node distance}
In order to study the distance between node $n$ and the root node $0_{1}$ in a $k$-tree it is natural to 
study the distance between node $n$ and all root nodes $0_{1}, \dots, 0_{k}$ simultaneously.
To do this we first introduce the notion $\text{dist}(u,K)$, which gives the distance between a node $u$ and a $k$-clique $K =\{w_{1}, \dots, w_{k}\}$ via $\text{dist}(u,K) := \min_{1 \le \ell \le k}\{\text{dist}(u,w_{\ell})\}$
(where $\text{dist}(x,y)$ denotes the distance between the nodes $x$ and $y$).
We introduce then the r.v. $\tilde{D}_{n}$, which counts the distance between node $n$ and the root-clique
$K_{0} = \{0_{1}, \dots, 0_{k}\}$ in a random $k$-tree of size $n$. 
Since the distance between node $n$ and an arbitrary root node in a $k$-tree is always either the same
as the distance between $n$ and the root-clique $K_{0}$ or one more,
there are always $\ell$ root nodes, with $1 \le \ell \le k$, which are at the same distance to node $n$ 
like the root-clique is, and $k-\ell$ root nodes with a distance one larger.
Due to symmetry it suffices to introduce the following $k$ different events $\mathcal{E}_{\ell}$, $1 \le \ell \le k$,
which describe the different situations that can occur:
\begin{equation*}
  \mathcal{E}_{\ell} := \big\{\text{dist}(n,0_{1}) = \cdots = \text{dist}(n,0_{\ell}) < \text{dist}(n,0_{\ell+1})
  = \cdots = \text{dist}(n,0_{k})\big\}.
\end{equation*}
Then the distribution of the r.v. $D_{n}$ we are interested in can be described as follows
(amongst the $\binom{k}{\ell}$ possible situations symmetric to event $\mathcal{E}_{\ell}$ one has to distinguish
whether node $0_{1}$ is at the same distance to $n$ or is one larger than the distance between $n$ and the root-clique):
\begin{equation}
   \label{eqns4_3}
   \mathbb{P}\{D_{n}=m\} = \sum_{\ell=1}^{k} \mathbb{P}\{\tilde{D}_{n} = m \wedge 
   \mathcal{E}_{\ell} \; \text{occurs}\} \binom{k-1}{\ell-1} 
   + \sum_{\ell=1}^{k} \mathbb{P}\{\tilde{D}_{n} = m-1 \wedge \mathcal{E}_{\ell} \; \text{occurs}\}
   \binom{k-1}{\ell}.
\end{equation}

When introducing the generating functions
\begin{align*}
   T_{\ell}(z,v) & := \sum_{n \ge 1} \sum_{m \ge 0} T_{n} \mathbb{P}\{\tilde{D}_{n} = m \; \wedge \;
   \mathcal{E}_{\ell} \; \text{occurs}\} \frac{z^{n-1}}{(n-1)!} v^{m}, \quad 1 \le \ell \le k, \\
   S_{\ell}(z,v) & := \sum_{n \ge 1} \sum_{m \ge 0} S_{n} \mathbb{P}\{\tilde{D}_{n}^{[S]} = m \; \wedge \;
   \mathcal{E}_{\ell} \; \text{occurs}\} \frac{z^{n-1}}{(n-1)!} v^{m}, \quad 1 \le \ell \le k,
\end{align*}
where $\tilde{D}_{n}^{[S]}$ denotes the corresponding r.v. for objects in the family $\mathcal{S}$,
one obtains by using the combinatorial decomposition of $k$-trees w.r.t. the root-clique 
(and after a study of the possibilities for the distance between node $1$ and node $n$ in objects of $\mathcal{S}$
leading to event $\mathcal{E}_{\ell}$)
given by \eqref{eqns2_3} the following system of equations, with $T(z)$ and $S(z)$ given in \eqref{eqns2_4}:
\begin{align*}
   T_{\ell}(z,v) & = \frac{S_{\ell}(z,v)}{(1-S(z))^{2}}, \quad 
   \frac{\partial}{\partial z} S_{\ell}(z,v) = (k-\ell) T(z)^{k-1} \big(T_{\ell}(z,v) + T_{\ell+1}(z,v)\big),
   \enspace 1 \le \ell \le k-1, \\
   T_{k}(z,v) & =  \frac{S_{k}(z,v)}{(1-S(z))^{2}}, \quad
   \frac{\partial}{\partial z} S_{k}(z,v) = k v T(z)^{k-1} T_{1}(z,v).
\end{align*}
This leads to the following system of linear differential equations for the functions $S_{\ell}(z,v)$:
\begin{equation*}
   \frac{\partial}{\partial z} S_{\ell}(z,v) = \frac{(k-\ell) \big(S_{\ell}(z,v)+S_{\ell+1}(z,v)\big)}{1-(k+1)z},
   \quad 1 \le \ell \le k-1, \quad \frac{\partial}{\partial z}S_{k}(z,v) = \frac{k v S_{1}(z,v)}{1-(k+1)z}.
\end{equation*}
Since it is possible to get from this system of differential equations a single differential equation for 
$S_{k}(z,v)$, which is of Euler type, it can be solved explicitly; thus all functions $T_{\ell}(z,v)$, 
$1 \le \ell \le k$, can also be given explicitly.

Since the generating function 
$T(z,v) := \sum_{n \ge 1} \sum_{m \ge 0} T_{n} \mathbb{P}\{D_{n} = m\} \frac{z^{n-1}}{(n-1)!} v^{m}$
is due to equation~\eqref{eqns4_3} completely determined by the functions $T_{\ell}(z,v)$, $1 \le \ell \le k$,
via
\begin{equation*}
   T(z,v) = \sum_{\ell=1}^{k} \left(\binom{k-1}{\ell-1} + \binom{k-1}{\ell} v\right) T_{\ell}(z,v),
\end{equation*}
it is possible to also get an explicit solution for $T(z,v)$.
One eventually obtains that
\begin{equation}
   T(z,v) = \sum_{j=1}^{k} \frac{B_{j}(v)}{(1-(k+1)z)^{\alpha_{j}(v)+\frac{2}{k+1}}},
   \label{eqns4_2}
\end{equation}
with $\alpha_{1}(v), \dots, \alpha_{k}(v)$ the roots of the equation
\begin{equation}
   \frac{k! v}{(k+1)^{k}} = \prod_{r=0}^{k-1}\big(\alpha - \frac{r}{k+1}\big),
   \label{eqns4_1}
\end{equation}
and where the functions $B_{j}(v)$, $1 \le j \le k$, can be given explicitly
(see \cite{Pan2004} for such considerations on a related problem);
for our purpose it is sufficient to state that all functions $B_{j}(v)$ are analytic in a neighbourhood of $v=1$.
By considerations as in \cite{Mah1992} one can show that for $v=1$ all roots $\alpha_{j}(1)$, $1 \le j \le k$,
of \eqref{eqns4_1} are simple. It is easily observed that $\frac{k}{k+1}$ is a root of \eqref{eqns4_1} when $v=1$; moreover, it is the root with largest real part. Let $\alpha_{1}(v)$ denote the root of \eqref{eqns4_1}, 
which satisfies $\alpha_{1}(1) = \frac{k}{k+1}$.
Then from \eqref{eqns4_2} we obtain the following asymptotic expansion of the moment generating function
of $D_{n}$:
\begin{equation*}
   \mathbb{E}\big(e^{D_{n} s}\big) = e^{(\alpha_{1}(e^{s})-\frac{k}{k+1}) \log n + 
   \log\big(\frac{\Gamma(1+\frac{1}{k+1}) B_{1}(e^{s})}{\Gamma(\alpha_{1}(e^{s})+\frac{2}{k+1})}\big)}
   \cdot \big(1+\mathcal{O}(n^{-\eta})\big), \quad \text{with an $\eta > 0$.}
\end{equation*}
An application of the quasi-power theorem of Hwang, see \cite{Hwa1998}, immediately shows Theorem~\ref{the7}.
Corollary~\ref{cor1} can be deduced from it easily.

\section{Conclusion}
We introduced a network model which is based on $k$-trees and which can either be described by a 
probabilistic growth rule using preferential attachment or combinatorially by considering increasing labellings
of the nodes and a linear ordering of the children of $k$-cliques. We gave a precise analysis of various parameters and could show that the distribution of the node-degrees follows asymptotically a power law, that the expected local clustering coefficient is high, and that the root-to-node distance of node $n$ is asymptotically Gaussian distributed
with expectation and variance of order $\log n$.
The approach we used is not restricted to the introduced model, but can also be easily applied to further evolving $k$-tree models; in particular the previously introduced uniform attachment model for $k$-trees can be treated in the same way. But also further $k$-tree models such as, e.g., ones with a ``saturation rule'', where at most up to $d$ children can be attached to a $k$-clique and where the probability that a new node is attached is proportional to the number of 
``free places'', can be introduced and analyzed (the special instance $d=1$ gives the so-called Apollonian networks).
One can even go a step further and introduce weighted ordered $k$-trees (with or without increasing labellings),
where each $k$-clique in the $k$-tree gets a weight depending on the number of attached children.
This, in analogy to simply generated tree families, see \cite{FlaSed2009}, leads then to ``simple families of $k$-trees'' (if unlabelled or arbitrarily labelled) or ``simple families of increasing $k$-trees'' (if increasingly labelled).
By choosing specific weights for increasingly labelled ordered $k$-trees all the before-mentioned evolution models for $k$-trees can be obtained; a complete characterization of possible $k$-tree evolution models as has been given in \cite{PanPro2007} for simple families of increasing trees is possible.

\nocite{*}
\bibliographystyle{alpha}
\bibliography{ktreesAofA10amsart}

\end{document}